\documentclass[11pt, reqno]{amsart}

\usepackage{amsmath,amsfonts,amsthm,amssymb,color,hyperref}

\usepackage{mathrsfs,soul}

\usepackage[T1]{fontenc}

\usepackage{graphicx}
\usepackage{subfigure} 

\makeatletter
     \def\section{\@startsection{section}{1}%
     \z@{.7\linespacing\@plus\linespacing}{.5\linespacing}%
     {\bfseries
     \centering
     }}
     \def\@secnumfont{\bfseries}
     \makeatother

\usepackage{graphicx}

\newcommand{\C}{\mathbb C}

\newcommand{\R}{\mathbb R}

\newcommand{\N}{\mathbb N}
\newcommand{\Q}{\mathbb Q}

\newcommand{\E}{\mathbb E}

\newcommand{\1}{{\bf 1}}
\newcommand{\e}{\varepsilon}

\newcommand{\HH}{\mathfrak H}

\newtheorem{theorem}{Theorem}[section]

\newtheorem{proposition}[theorem]{Proposition}

\theoremstyle{definition}
\newtheorem{definition}[theorem]{Definition}
\theoremstyle{remark}
\newtheorem{remark}{Remark}
\numberwithin{equation}{section}
\begin{document}

\title[Oscillatory Breuer-Major theorem]{Oscillatory Breuer-Major theorem with application to the random corrector problem}

 \date{\today}

\author{David Nualart \and Guangqu Zheng}

\address{David Nualart: University of Kansas, Mathematics department, Snow Hall, 1460 Jayhawk blvd, Lawrence, KS 66045-7594, United States} \email{nualart@ku.edu} \thanks{D. Nualart is supported by NSF Grant DMS 1811181.}

\address{Guangqu Zheng: University of Kansas, Mathematics department, Snow Hall, 1460 Jayhawk blvd, Lawrence, KS 66045-7594, United States} \email{zhengguangqu@gmail.com}

\begin{abstract}
In this paper, we present an oscillatory version  of the celebrated Breuer-Major theorem that is motivated by the random corrector problem. As an application, we are able to prove new results concerning the Gaussian fluctuation of  the random corrector.   We also provide a variant of this theorem involving homogeneous measures.

\end{abstract}
 
 \maketitle

  \section{Introduction and main results}
  
  Our work is motivated by the following random homogenization problem. Consider a  one-dimensional
  equation with highly {\it oscillatory} coefficients of the form
 \begin{align}\label{EE}
 \begin{cases}
 \displaystyle{ -\frac{d}{dx}\left( \,\, a(x/\varepsilon, \omega) \frac{d}{dx} u_\e(x, \omega) \, \right) } = f \in L^1([0,1], dx)  \\
 u_\e(0, \omega) = 0 \,\,,\quad u_\e(1, \omega) = b\in\R,
 \end{cases}
 \end{align}
 where $\varepsilon \in (0,1]$.
In the literature (see \emph{e.g.} \cite{BGMP08, BG15, GB12, JKO94}), the random potential $a$ is often assumed to be ergodic, uniformly elliptic (\emph{i.e.} positive and bounded with bounded inverse). Notice that, under the following hypothesis:
 \begin{center}    \quad    For all $\varepsilon \in (0,1]$,  ${\displaystyle \int_0^{1/ \varepsilon}  \frac 1{|a(x)|} dx <\infty}$  and  ${\displaystyle \int_0^{1/ \varepsilon}  \frac 1{a(x)} dx  \not =0}$  almost surely,    \hfill  {\bf (H)}
 \end{center}
we can solve \eqref{EE} explicitly:
\begin{align}\label{sole}
 u_\e(x, \omega) = c_\e(\omega) \int_0^x \frac{1}{a(y/\varepsilon,\omega)} \, dy - \int_0^x \frac{F(y)}{a(y/\varepsilon, \omega)} \, dy,
\end{align} 
 where
 $F(x) := \int_0^x f(y) \, dy$
 is the antiderivative of $f$ vanishing at zero
 and 
 $$c_\e(\omega) : = \left( \, b + \int_0^1 \frac{F(y)}{a(y/\varepsilon, \omega)} \, dy \,\right) \left(  \int_0^1 \frac{1}{a(y/\varepsilon, \omega)} \, dy  \right)^{-1}.$$

 {\bf Throughout this note}, we    assume  that $a$  satisfies {\bf (H)} and    has the following  form 
 \begin{align}\label{aform}
  a(x)  =  \left( \frac{1}{a^\ast}  + \Phi(W_x)  \right)^{-1}\,,
 \end{align}
 where 
 
  (i)
   $\{ W_x , x\in\R \}$ is a centered {\it stationary Gaussian process} with a   correlation given by $\rho(x-y) = \E\big[ W_x W_y \big]$, and we assume that $\rho$ is {\bf continuous} with $\rho(0)=1$; 
   
    (ii)  $\Phi\in L^2(\R, e^{-x^2/2}dx )$  has the following orthogonal expansion 
 \begin{align} \label{PHI}
 \Phi(x) = \sum_{q\geq m} c_q H_q \,,
 \end{align}
 with $H_q(x)= (-1)^q e^{x^2/2} \frac{d^q}{dx^q} e^{-x^2/2}$ denoting the $q$th Hermite polynomial.  Here $c_m\neq 0$ and $m\geq 1$ is called the {\it Hermite rank} of $\Phi$.  
  
The quantity  $a^\ast :=1/ \E [ 1/a(0)]$ is known as the \emph{harmonic mean} or {\it effective diffusion coefficient} of the random potential, see \cite{JKO94,PV81}.

\begin{remark}Assuming the structure $a(x)^{-1} = \Phi(W_x) + (a^\ast)^{-1} $,  our hypothesis   {\bf (H)}  holds  provided $\int_0^{1/ \varepsilon}  \left(  (a^\ast)^{-1}  + \Phi(W_x) \right)dx  \not =0$  almost surely, for all $\varepsilon \in (0,1]$.  Note that  the local integrability of $a(x)^{-1}$ follows immediately from its  structure: Indeed, for any $\e > 0$,
\[
\E \int_0^{1/\e} \frac{1}{\vert a(x) \vert} dx \leq \frac{1}{\vert a^\ast \vert \e} + \int_0^{1/\e} \E\big[ \vert \Phi(W_x) \vert \big]dx = \frac{1}{\vert a^\ast \vert \e} +   \frac{1}{ \e}\E\big[\vert \Phi(W_1) \vert \big]  < +\infty,
\]
which implies that $ \int_0^{1/\e}  \vert a(x) \vert^{-1} dx $ is almost surely finite.  It is clear that our hypothesis   {\bf (H)} holds in presence of uniform ellipticity of $a$ and the latter is  equivalent to the boundedness of $\Phi$; as one can see from page 276-277 in \cite{Atef}, one can easily construct bounded measurable function $\Phi$ with given Hermite rank.  Note that if $\Phi(x)= \vert x\vert - \sqrt{2/\pi}$ (this is unbounded with Hermite rank $2$) and $a^\ast = \sqrt{\pi/2}$, then $a(x) = \vert W_x \vert^{-1}$
satisfies the assumption {\bf (H)} but not the  uniform ellipticity.
\end{remark}

 Under some mild assumptions on $\rho$ and $\Phi$, we can derive the following result concerning the asymptotic behavior of $u_\e$ as well as the associated fluctuation.

\begin{theorem}\label{KU}  Let the above notation prevail.  We assume that  $a$ satisfies {\bf (H)} and has the form \eqref{aform} such that $\Phi$, given as in \eqref{PHI}, has Hermite rank $m\geq 1$ and   the  correlation function $\rho$ of  the centered stationary Gaussian process $\{W_x,x\in\R\}$ belongs to $L^m(\R, dx)\cap C(\R)$ with $\rho(0)=1$. Then the following statements hold true:
\begin{itemize}
\item[(1)] 

For every $x\in[0,1]$,  $u_\e(x)$ converges  in probability to  $\bar{u}(x)$, as $\e\downarrow 0$, where $\bar{u}(x)$ solves the following $($deterministic$)$ homogenized equation
 \begin{align}\label{EEh}
 \begin{cases}
 \displaystyle{ -\frac{d}{dx}\left( \,\, a^\ast \frac{d}{dx} \bar{u}(x) \, \right) } = f  \\
 \bar{u}(0) = 0 \,\,,\qquad \bar{u}(1) = b .
 \end{cases}
 \end{align}

\item[(2)] For every $x\in(0,1)$, with $\mu^2:= \sum_{q\geq m} c_q^2q! \int_{\R} \rho(t)^qdt\in[0,\infty)$,
\begin{align}\label{goal1}
\frac{u_{\e}(x) - \bar{u}(x)}{\sqrt{\e}} \xrightarrow[\rm law]{\e\downarrow 0} N\left(0, \,\, \mu^2\int_0^1 F(x,y)^2 dy\right)\,.
\end{align}
 Moreover, if in addition $\Phi\in L^p(\R, e^{-x^2/2}dx)$ for some $p>2$, then 
\begin{align}\label{goal2}
\left\{ \frac{u_{\e}(x) - \bar{u}(x)}{\sqrt{\e} }\,, \, x\in[0,1] \right\} \xrightarrow[ \rm law]{\e\downarrow 0 }  \left\{ \mu \int_{0}^1 F(x,y) dA_y\,, x\in [0,1] \right\},
\end{align}
where  the above weak convergence takes place in $C([0,1])$, 
\[
F(x,y) := \big( c^\ast - F(y) \big)\1_{[0,x]}(y) + x \big( F(y) - c^\ast\big)
\]
 for $x,y\in[0,1]$ and $\{A_y, y\in[0,1]\}$ is a standard Brownian motion. Here $c^\ast:=ba^\ast + \int_0^1 F(z)dz$.

\end{itemize}

\end{theorem}
The difference $u_\e - \bar{u}$ is known as the {\it random corrector} in the homogenization theory, see \cite{BGMP08} and references therein.  Our Theorem \ref{KU} complements   findings   in the literature, see the following Remark \ref{rem21}: Points (i)-(iii) sketch some relevant history and points (iv)-(v) summarize the novelty of our results.
  
   \begin{remark}\label{rem21} (i) The authors of \cite{BP99} considered the  short-range case where the random potential  $\{a(x), x\in\R\}  $  satisfies certain (strong) mixing conditions:  With the above notation,  mixing conditions and uniform ellipticity  in \cite{BP99}  imply  that $\E[ \Phi(W_x) \Phi(W_y) ] $  is bounded by  $\text{constant} \cdot \vert x-y\vert^{-\alpha}$ for some $\alpha>1$. Since  the correlation function of $\Phi(W_x)$    is also bounded by $\| \Phi\|_\infty^2$, it is integrable, which guarantees that the random corrector $u_\e - \bar{u}$ is of order $\sqrt{\e}$; properly scaled, the random corrector  converges to a Wiener integral with respect to Brownian motion; see also Theorem 2.6 in \cite{BGMP08}.

 (ii) In \cite{BGMP08}, the result has been extended to a large family of random potential with long-range correlation (\emph{i.e.} $\rho(\tau)\sim \text{constant}\cdot \tau^{-\alpha}$ for some $\alpha\in(0,1)$): It was shown that when the Hermite rank of $\Phi$ is one,  the corrector's amplitude is of order $\e^{\alpha/2}$ and after properly scaled, the random corrector converges in law to a stochastic integral with respect to the fractional Brownian motion with Hurst parameter $(2-\alpha)/2$; see also Theorem 2.3 in \cite{GB12}.
 
 (iii) Following \cite{BGMP08}, the authors of \cite{GB12} studied the random corrector problem for the case where the Hermite rank of $\Phi$ is two and $\rho(\tau)\sim \text{constant}\cdot \vert\tau\vert^{-\alpha}$ as $\tau\to\infty$, with $\alpha\in(0,1/2)$. They established that the corrector's amplitude is of order $\e^{\alpha}$ and the random corrector, after proper rescaling, converges in law to a stochastic integral with respect to the Rosenblatt process; see \cite[Theorem 2.2]{GB12}. In the end of the paper \cite{GB12}, the authors conjectured that when the Hermite rank of $\Phi$ is three or higher, the properly rescaled corrector is expected to converge in law to some stochastic integral with respect to the so-called Hermite process and this is confirmed in the work \cite{Atef}.

 (iv) Note that all the references mentioned in  (i)-(iii) assume that $a$ is stationary ergodic such that $0< c_1 \leq a(x) \leq c_2  $ almost surely for some   numerical constants $c_1, c_2$ (so $\Phi$ is bounded), while  we do not assume the uniform boundedness of $\Phi$. Instead, we only assume   hypothesis {\bf (H)} and $\Phi\in L^{p}(\R, e^{-x^2/2}dx)$ for some $p>2$.    
 In our framework, the correlation function $\rho$   belongs to $ L^m(\R,dx)$, with $m\geq 1$ being the Hermite rank of $\Phi$, which ensures that  the correlation function of  $\{\Phi(W_x),x\in\R\}$ is integrable. So similar to  \cite{BP99}, we are in the short-range setting and we establish that the corrector's amplitude is of order $\sqrt{\e}$ and properly rescaled corrector converges in law to a Gaussian process.

(v) For the functional convergence \eqref{goal2}, we impose the condition $p>2$    in order to have moment estimates of order $p$ that imply tightness. These moment estimates are derived  using  Meyer's inequalities.   The aforementioned example $\Phi(x) = \vert x\vert  -\sqrt{2/\pi}$ belongs to $L^p(\R, e^{-x^2/2}dx)$ for any $p\geq 2$. Our proof of Theorem \ref{KU} uses techniques from Malliavin calculus and Gaussian analysis, which might be helpful for other more complicated problems in  random PDEs.

 \end{remark}

Our  Theorem \ref{KU} is a special case of the following more general result. 
We denote by $\mathcal{B}_b$ the collection of bounded  closed sets in $\R^d$. For any $R\ge 0$ we put  $B_R:=\{ x\in\R^d\,:\, \| x\| \leq R \}$.  Also, {\it f.d.d.} means convergence  of the finite-dimensional distributions of a given family of random variables depending on a parameter $R$, which tends to $+\infty$.

  \begin{theorem} \label{Lawrence} Let $\{ W_x, x\in \R^d\}$ be a centered  Gaussian stationary process with  continuous  covariance
$
\rho(x-y) := \E[ W_x W_y]
$ such that $\rho(0)=1$ and $\rho\in L^m(\R^d, dx)$. Let $\Phi$ be given as in \eqref{PHI} with Hermite rank $m\geq 1$.
Then, with $h\in C(\R^d)$,  we have
\begin{align}\label{CLT1}
\left\{ R^{d/2}\int_{B} \Phi(W_{xR}) h(x) \, dx \right\}_{ B\in \mathcal{B}_b}  \xrightarrow[\rm f.d.d.]{R\to+\infty}  \left\{ \sigma\int_{B}  h(x) dZ_x \right\}_{ B\in \mathcal{B}_b},
\end{align}
where $Z$ denotes the standard Gaussian white noise on $\R^d$ and  
\[
\sigma^2=   \sum_{q=m} ^\infty q! c_q^2 \int_{\R^d} \rho(z)^q \, dz \in[0,+\infty)\,.
\]

If in addition $\Phi\in L^p(\R, e^{-x^2/2}dx)$ for some $p>2$. 
Then,  the following functional central limit theorems hold true: 

{\rm (1)}  With  $p>2d$ and any finite $\ell > 0$,
\begin{align}\label{FCLT20}
\left\{ R^{d/2} \int_{[0,\pmb{z}]}  \Phi(W_{xR}) h(x) dx \right\}_{\pmb{z}\in [0,\ell]^d} \xrightarrow[R\to+\infty]{\rm law}    \left\{ \sigma \int_{[0,\pmb{z}]}  h(x) dZ_x\right\}_{\pmb{z}\in   [0,\ell]^d} ,
\end{align}
where     the above weak convergence holds on the space $C\big([0,\ell]^d\big)$ and  $[0,\pmb{z}]= \prod_{j=1}^d [0,z_j]$ given $\pmb{z} = (z_1, \ldots, z_d)\in [0,\ell]^d$; 
  
{\rm (2) }   
\begin{align}
\left\{ R^{ d/2} \int_{B_{t} } \Phi(W_{xR}) h(x) dx \right\}_{t\ge 0} \xrightarrow[R\to\infty]{\rm law}   \left\{ \sigma\int_{B_t}  h(x) dZ_x\right\}_{t\geq 0}, \label{FCLT1}
\end{align}
where  the above weak convergence takes place on $C(\R_+)$.

\end{theorem}
 
 Roughly speaking, the random corrector $u_\e(x) - \bar{u}(x)$ from Theorem \ref{KU} can be written as a sum of an oscillatory integral and a negligible  term so that an easy application of  Theorem \ref{Lawrence} gives us Theorem \ref{KU}, see Section \ref{SEC2} for more details.  We will   proceed the proof of \eqref{CLT1} by following the usual arguments for the chaotic central limit theorem (see \emph{e.g.} \cite{HN05,bluebook}), while the functional central limit theorem in \eqref{FCLT1} is established with the help of Malliavin calculus techniques, notably  Meyer's inequality (see \cite{CNN18,NN18}).  
 
 \begin{remark} (i) Theorem \ref{Lawrence} is a generalization of the celebrated Breuer-Major theorem \cite{BM83} that corresponds to the case where $h=1$, see also \cite{CNN18, NN18}.   The integral on the left-side of \eqref{CLT1} is known as an oscillatory random  integral, 
  so we call our result an oscillatory Breuer-Major theorem and this explains our title. 
 
 (ii) The functional limit theorem described in \eqref{FCLT20} is new and  the limit is a $d$-parameter Gaussian process  with covariance given by 
$$
 \sigma^2\int_{[0,\pmb{z}] \cap [0,\pmb{y}]} h^2(x) dx\,,
$$
while  the limit in \eqref{FCLT1} is a Gaussian martingale with quadratic variation given by $$ t\in\R_+ \longmapsto \sigma^2\int _{B_t}  h^2(x) dx\,.$$
 
 \end{remark}
 \medskip
  
  Our approach is quite flexible and  we can provide another variant of Breuer-Major's  theorem that involves an  {\it homogeneous measure}. Let us first recall the definition of homogeneous measure (see \emph{e.g.} \cite{EM18}).

\begin{definition}  Given $\alpha\in\R\setminus\{0\}$,  a measure $\nu$ on $\R^d$ is said to be $\alpha$-{\it homogeneous} if 
\begin{center}
$\nu(sA) = s^{\alpha} \nu(A)$, for any $s > 0$ and $A\subset \R^d$ Borel measurable,
\end{center}
where $sA := \big\{ x\in\R^d\,: s^{-1}x\in A\big\}$. For example, $\mu(dx)= \vert x\vert^{-\beta}dx$ defines a $(d-\beta)$-homogeneous measure on $\R^d$ for any $\beta\neq d$. Note that for general $h\in C(\R^d)$, the measure $\gamma(dx) = h(x)dx$ is not necessarily homogeneous. 
\end{definition}

\begin{theorem}\label{Kansas} Fix $\alpha\in(0,\infty)$ and consider an $\alpha$-homogeneous measure $\nu$ on $\R^d$ such that $0 < \nu(B_1) <\infty$.  Let $\Phi$ be given as in \eqref{PHI} with Hermite rank $m\geq 1$ and let $\{ W_x, x\in \R^d\}$ be a centered  Gaussian stationary process with  continuous  covariance
$
\rho(x-y) := \E[ W_x W_y]
$ such that $\rho(0)=1$ and $\rho\in L^m(\R^d, d\nu)$. Then
\begin{align}\label{CLT2}
 \left\{  R^{\alpha/2} \int_{B} \Phi(W_{xR}) \nu(dx) \right\}_{ B \in \mathcal{B}_b}
  \xrightarrow[\rm f.d.d.]{R\to+\infty}   \left\{ \sigma_\nu \, \mathcal{Z}(B) \right\}_{B \in \mathcal{B}_b}
\end{align}
where  $\mathcal{Z}$ stands for the Gaussian random measure with intensity   $\nu$ on $\R^d$  and 
$$\sigma_\nu^2 :=\sum_{q\geq m} c_q^2 q! \int_{\R^d} \rho(z)^q\, \nu(dz)\in [0,+\infty).$$ Moreover,    if   additionally  $\Phi\in L^p(\R, e^{-x^2/2}dx)$ for some $p>2$ and $\alpha p > 2$, then we have the following functional central limit theorem:
  \[
\left\{ R^{ \alpha/2} \int_{B_{t} } \Phi(W_{xR}) \nu(dx) \right\}_{t\ge 0} \xrightarrow[R\to+\infty]{\rm law}  \Big\{ \sigma_\nu \, \mathcal{Z}(B_t) \Big\}_{t\ge 0}.
\]
\end{theorem}

 One can refer to the book \cite{Nualart06} for any unexplained notation and definition.
 We would like to point out that if  the function $\Phi$ is a finite sum of Hermite polynomials, then, the Stein-Malliavin approach implies
  that, in the framework Theorem \ref{Kansas},  the convergence of the one-dimensional distributions hold  in  the total variation distance
  (see for instance, the monograph  \cite{bluebook}).
  
The rest of this article  consists of three more sections: Section \ref{PRE} is devoted to some    preliminary material. In Section \ref{SEC2}, we present the proof of Theorem \ref{Lawrence} and then as anticipated, we demonstrate how Theorem \ref{Lawrence} implies Theorem \ref{KU}.   We will sketch the proof of Theorem \ref{Kansas} in Section \ref{SEC4}.  

Note that all random objects in this note are assumed to be defined on a common probability space $(\Omega, \mathscr{F}, \mathbb{P})$ and we will use $C$ to denote a generic constant that is immaterial to our estimates and it may vary from line to line.

  \section{Preliminaries}\label{PRE}

  Recall that $\{W_x,x\in\R^d\}$ is a centered  stationary Gaussian process such that it has  a    continuous   covariance function $\rho$.  The continuity of $\rho$ is equivalent to the $L^2(\Omega)$-continuity of process $W$. In what follows, we first build the isonormal framework  for later Gaussian analysis.  Note that the Gaussian Hilbert space generated by $W$ is the same as the one generated by $\{W_x,x\in\Q^d\}$ due to the $L^2$ continuity, so the resulting Gaussian Hilbert space is a real separable Hilbert space. By a standard fact in real analysis, it is isometric to $L^2([0,1], dt)=:\HH$ and we denote this isometry by $X$. By isometry, there exists a sequence $\{ e_x, x\in\Q^d\}\subset \HH$ such that
  \[
  X(e_x) = W_x \quad\text{for any $x\in\Q^d$.}
  \]
  By continuity again, the above equality extends to every $x\in\R^d$. It is clear that $\{X(h), h\in\HH \}$ is  an isonormal Gaussian process over the  real separable Hilbert space $\HH$. By construction,  $e_x\in \HH$ has unit norm and $\langle e_x, e_y \rangle_{\HH} = \rho(x-y)$ for any $x,y\in\R^d$. Note that $x\in\R^d\longmapsto e_x\in \HH$ is a continuous map and this can save us away from measurability issues.

 In what follows, we introduce some standard notation from Malliavin calculus; see the basic references  \cite{bluebook,Nualart06, Eulalia} for more details.   For a smooth and cylindrical random variable $F = f\big(X(h_1), \ldots, X(h_n)\big)$ with $h_i\in\HH$ and $f\in C^\infty_b(\R^n)$, we define its Malliavin derivative as the $\HH$-valued random variable given by 
  \[
  DF =  \sum_{i=1}^n \frac{\partial}{\partial x_i} f\big(X(h_1), \ldots, X(h_n)\big) h_i.
  \]
By iteration, we can define the $k$th Malliavin derivative of $F$ as an element in $L^2(\Omega; \HH^{\otimes k})$.  Here $\HH^{\otimes k}$ denotes the $k$th tensor product of $\HH$ and we denote by $\HH^{\odot k}$ the space of symmetric tensors in $\HH^{\otimes k}$. For any  $k\in\N$ and $p\in[ 1,\infty)$, we define the Sobolev space $\mathbb{D}^{k,p}$
as the
closure of the space of smooth and cylindrical random variables with respect to the norm $\| \cdot\| _{k,p}$ defined by
\[
\| F\| _{k,p}^p = \E\big( | F |^p \big) + \sum_{i=1}^k \E\big( \| D^i F \| _{\HH^{\otimes i}}^p \big) \,.
\]
The divergence operator $\delta$ is defined as the adjoint of the derivative operator $D$. An element $u\in L^2(\Omega; \HH)$
  belongs to the domain of $\delta$,  denoted by $\text{dom}(\delta)$  if there is a constant $c_u$ that only depends on $u$ such that 
  \[
  \big\vert \E\big[ \langle DF, u \rangle_\HH \big] \big\vert \leq c_u  \sqrt{\E [ F^2 ]} \quad\text{for any $F\in\mathbb{D}^{1,2}$.}
  \]
For $u\in \text{dom}(\delta)$, the existence of $\delta(u)$ is guaranteed by the Riesz representation theorem and it satisfies the following duality relation
\[
\E\big[ \langle DF, u \rangle_\HH \big] = \E\big[ F \delta(u) \big] \quad\text{for any $F\in\mathbb{D}^{1,2}$.}
\]
Similarly, we can define the iterated divergence $\delta^k$: For $u\in\text{dom}(\delta^k)\subset L^2(\Omega; \HH^{\otimes k})$, $\delta^k(u)$ is characterized by the following duality relation
  \[
 \E\big[ \langle D^kF, u \rangle_{\HH^{\otimes k}} \big] = \E\big[ F \delta^k(u) \big] \quad\text{for any $F\in\mathbb{D}^{k,2}$.}
  \]
The well-known Wiener-It\^o chaos decomposition states that  any $F\in L^2(\Omega, \sigma\{W\}, \mathbb{P})$ admits the following expression 
\begin{align}\label{chaos}
F = \E[F] + \sum_{p\geq1} \delta^p(f_p),
\end{align}
with $f_p\in\HH^{\odot p}$ uniquely determined by $F$;  $ \delta^p(f_p)$ is also called the $p$th multiple integral with kernel $f_p$.    Note that  given any   unit vector $e\in\HH$, we have $H_p(X(e)) = \delta^p(e^{\otimes p})$. We call $\mathbb{C}_p$, the closed linear subspace of $L^2(\Omega)$ generated by  $\big\{ H_p(X(e)) \,:\, e\in\HH$ and $\| e\| _\HH=1\big\}$, the $p$th   Wiener chaos associated with the isonormal Gaussian process $X$ and we write $J_p$ for the projection operator onto $\mathbb{C}_p$.
Then we define Ornstein-Uhlenbeck semigroup $(P_t, t\in\R_+)$ and its generator $L$ by putting
\[
P_t = \sum_{p\geq 0} e^{-pt} J_p \quad \text{and} \quad L = \sum_{p\geq 1} -p J_p\,,
\]
and we write $L^{-1}$ for the pseudo-inverse of $L$, that is, 
\[
L^{-1}F =  -\sum_{p\geq 1}  \frac{1}{p} J_pF \quad\text{for any centered $F\in L^2(\Omega, \sigma\{W\}, \mathbb{P})$.}
\]
Note that these operators enjoy the following nice relation:  $F = - \delta D L^{-1} F$ for any centered $F\in L^2(\Omega, \sigma\{W\}, \mathbb{P})$. Now let us record an important consequence of this relation. Let $\Phi$ be given as in \eqref{PHI} and have Hermite rank $m\geq 1$. We define the shifted function
\[
\Phi_m(x)  = \sum_{q\geq m} c_q H_{q-m}(x) \,,
\]
which satisfies the following properties:

\noindent{(\pmb{A})} $\Phi_m(W_x) = \Phi_m\big(X(e_x)\big)   \in \mathbb{D}^{m,2}$ and $\Phi(W_x) = \delta^m\big( \Phi_m(W_x) e_x^{\otimes m} \big)$ for any $x\in\R^d$; 

\noindent{(\pmb{B})} $\Phi_m(W_x)e_x^{\otimes m} = (-DL^{-1})^m \Phi(W_x)$ and   applying Meyer's inequality, we have for every $k\in\{0,1, \ldots, m\}$, $x\in\R^d$ and  $p>1$,
\begin{align}\label{IMP}
 \big\| D^k \big(\Phi_m(W_x) \big) \big\| _{L^p(\Omega; \HH^{\otimes k})} \leq C \| \Phi(W_x) \|_{L^p(\Omega)}\,.
\end{align}
   This inequality  is a consequence of   Lemmas 2.1, Lemma 2.2 in \cite{NN18} (see also   \cite[(2.7)]{CNN18}).

 Let $\{\e_i, i\in \mathbb{N}\}$ be an orthonormal basis of $\HH$. For $f\in\HH^{\odot p}$ and $g\in\HH^{\odot q}$ ($p,q\in\N$), we define the $r$-contraction as the element in $\HH^{\otimes p+q-2r}$  ($r\in\{0,\ldots, p\wedge q\}$) given by
\[
f\otimes_r g =  \sum_{i_1, \ldots, i_r\in \mathbb{N}} \big\langle f , \e_{i_1} \otimes  \e_{i_2} \otimes \cdots  \otimes  \e_{i_r}  \big\rangle_{\HH^{\otimes r}}  \big\langle g , \e_{i_1} \otimes  \e_{i_2} \otimes \cdots  \otimes  \e_{i_r}  \big\rangle_{\HH^{\otimes r}} \,.
\]
  In particular, $f\otimes_0 g  = f\otimes g$ and if $p=q$, $f\otimes_p g  = \langle f, g\rangle_{\HH^{\otimes p}}$. 
  
  \medskip
  
  In the end of this section, we present a multivariate version of the chaotic central limit theorem \cite{HN05} that we borrow from \cite[Theorem 2.1]{CNN18}.

\begin{proposition}\label{CCLT} Fix an integer $n\geq 1$ and consider a family $\big\{ G_R, R> 0 \big\}$ of  random vectors in $\R^n$ such that each component of $G_R=(G_{R, 1}, \ldots, G_{R,n})$ belongs to $ L^2(\Omega, \sigma\{W\}, \mathbb{P})$ and has the following chaos expansion
\[
G_{R,j} = \sum_{q\geq 1} \delta^q\big(  g_{q, j,R} \big)\quad \text{with $g_{q, j,R}\in \HH^{\odot q}$ deterministic. }
\]
Suppose the following conditions {\rm (a)-(d)} hold:
\begin{itemize}
\item[(a)] For each $i,j\in\{1, \ldots, n\}$ and for every $q\geq 1$, $q! \langle g_{q, i,R}, g_{q, j,R}\rangle_{\HH^{\otimes q}}$ converges to some $\sigma_{i,j,q}\in\R$, as $R\to+\infty.$

\item[(b)] For each $i\in\{1, \ldots, n\}$, $\sum_{q\geq 1} \sigma_{i,i,q} < +\infty$.

\item[(c)] For each  $i\in\{1, \ldots, n\}$,  $q\geq 2$ and $r\in\{1, \ldots, q-1\}$, we have that, as $R\to+\infty$,  $\big\| g_{q, i,R}\otimes_r g_{q, i,R} \big\| _{\HH^{\otimes 2q-2r}}$ converges to zero.

\item[(d)] For each  $i\in\{1, \ldots, n\}$, $\lim_{N\to+\infty} \sup_{R>0} \sum_{q\geq N+1} q! \| g_{q, i,R} \|^2 _{\HH^{\otimes q}} = 0$.

\end{itemize}
Then $G_R$ converges in law to $N(0, \Sigma)$ as $R\to+\infty$, where $\Sigma = \big( \sigma_{i,j} \big)_{i,j=1}^n$ is given by  $\sigma_{i,j} = \sum_{q\geq 1} \sigma_{i,j,q}$.
  \end{proposition}
 
 The above proposition is essentially a consequence of the Fourth Moment Theorems due to Nualart, Peccati and Tudor  (see \cite{FMT, PT05}): In 2005, Nualart and Peccati discovered that for $\{F_n,n\geq 1\}\subset \C_p$ ($p\geq 2$), if $\E[F_n^2]\to 1$, then the asymptotic normality of this sequence is equivalent to  $\E[F_n^4]\to3$. Soon later, Peccati and Tudor provided a multidimensional extension, which asserts that for a sequence of random vectors  $G_n=(G_{1,n}, \ldots, G_{d,n})$ with covariance matrix convergent to some covariance matrix $C$, if for each $j$, $p_j\geq 1$,   $\{ G_{j,n}: n\geq 1 \}\subset\C_{p_j}$, then the joint convergence ($G_n$ converges in law to $N(0, C)$) is equivalent to the marginal convergence  ($G_{j,n}$ converges in law to $N(0, C_{jj})$ for each $j$).  The latter boils down to checking the fourth moment condition. For example, in the setting of Proposition \ref{CCLT}, let us look at the convergence of $G_{R,1}$:   conditions (b) and (d) ensure that it suffices to consider finite many chaoses, conditions (a) and (b) guarantee the convergence of the covariance matrix of the random vectors formed by these finitely many chaoses. In view of the product formula for multiple integrals, verifying the fourth moment condition would lead to the computation involving the contractions, where we need condition (c) for this to work; see \cite{CNN18} for a proof and we refer the interested readers to  the monograph \cite{bluebook} for a comprehensive introduction to this line of research.

    \section{Proof of Theorem \ref{Lawrence}   and  Theorem \ref{KU}  }\label{SEC2}
    
    In this section, we first prove the  convergence of finite-dimensional distributions in the framework of Theorem \ref{Lawrence}. Next, we will establish the tightness property under the additional assumption that $\Phi\in L^p(\R, e^{-x^2/2}dx )$ for some $p>2$, which is needed
    to establish \eqref{FCLT20} and \eqref{FCLT1}. These two steps will conclude the proof of Theorem \ref{Lawrence}, and in the end of this section, we demonstrate how one can derive Theorem \ref{KU} from Theorem \ref{Lawrence}.

  \subsection{Convergence of finite-dimensional distributions}\label{FDDcvg} For each $R > 0$ and $B\in \mathcal{B}_b$, we put
  \begin{align}\label{GRT}
  G_R(B) = R^{d/2}\int_{B} \Phi(W_{xR}) h(x) \, dx\,.
  \end{align}
  Then,  it is enough to consider bounded Borel sets $B_i \in \mathcal{B}_b$, $i=1,\dots, n$,  and establish the following limit result 
  \begin{align} \label{FDD}
  \big( G_R(B_1), \ldots, G_R(B_n) \big)  \xrightarrow[\rm law]{R\to+\infty} N(0, \Sigma) \,,
  \end{align}
  where $\Sigma = \big( \sigma_{i,j}\big)_{i,j=1}^n$ is defined by $$\sigma_{i,j} = \sigma^2 \int_{B_i \cap B_j } h(x)^2\, dx \,.$$
For $j\in\{1,\ldots, n\}$, we can rewrite $G_R(B_j)$ using the Hermite expansion \eqref{PHI} as follows:
\begin{align*}
G_R(B_j) &=  R^{d/2}\int_{B_j} \sum_{q\geq m} c_q H_q(W_{xR}) h(x) \, dx = R^{d/2}\int_{B_j} \sum_{q\geq m} \delta^q(c_q e_{xR}^{\otimes q}) h(x) \, dx \\
& = \sum_{q\geq m}  \delta^q\left(  c_q R^{d/2} \int_{B_j}e_{xR}^{\otimes q} h(x) \, dx \right) =:\sum_{q\geq m}  \delta^q\left(  g_{q,j,R} \right) \,,
\end{align*}
 where
 \[
 g_{q,j,R} = c_q R^{d/2} \int_{B_j}e_{xR}^{\otimes q} h(x) \, dx.
 \]
   
(a) For any $i,j\in\{1, \ldots, n\}$, we have 
\begin{align*}
q! \langle g_{q, i,R}, g_{q, j,R}\rangle_{\HH^{\otimes q}} &= q! c_q^2 R^d \int_{ B_i \cap B_j } \rho(xR-yR)^q h(x)h(y) \, dxdy \\
& = q! c_q^2 R^{-d}\int_{ RB_i \cap RB_j } \rho(x-y)^q h(x/R)h(y/R) \, dxdy \\
&=  q! c_q^2 R^{-d}\int_{  \{x\in RB_i , x-z\in RB_j\} } \rho(z)^q h(x/R)h((x-z)/R) \, dxdz.
  \end{align*}
Making the change of variables $x/R= y$ yields
\begin{align*}
q! \langle g_{q, i,R}, g_{q, j,R}\rangle_{\HH^{\otimes q}}  &  =   c_q^2 q!  \int_{  \{ y\in B_i, y-zR^{-1} \in B_j \}} \rho(z)^q  h\big( y-zR^{-1}\big) h(y) \, dydz.
\end{align*}
Taking into account  that  $h$ is continuous and $B_j$ is closed,    we   deduce from the dominated convergence theorem that 
\[
q! \langle g_{q, i,R}, g_{q, j,R}\rangle_{\HH^{\otimes q}} \xrightarrow{R\to+\infty} c_q^2 q! \left(\int_{\R^d} \rho(z)^q dz\right)   \int_{B_i \cap B_j}   h(y)^2 \, dy =: \sigma_{i,j,q}\,.
  \]
  
  (b) For each $i\in\{1, \ldots, n\}$, 
  $$
  \sum_{q\geq m} \sigma_{i,i,q} = \left(\int_{B_i}   h(y)^2 \, dy\right)  \sum_{q\geq m}  c_q^2 q! \left(\int_{\R^d} \rho(z)^q dz\right) = \sigma^2\int_{B_i}   h(y)^2 \, dy  \,.
  $$
  Note that the quantity $\sigma^2$ as defined in the statement of Theorem \ref{Lawrence} is finite, because $\int_{\R^d} \rho(z)^q dz$ is bounded by $\int_{\R^d} \vert \rho(z)\vert^mdz$ and $\sum_{q\geq m}c_q^2 q! <+\infty$. So we just verified the condition (b).
  
  \medskip
  
  (c) For each  $i\in\{1, \ldots, n\}$,  $q\geq 2$ and $r\in\{1, \ldots, q-1\}$, we have,  
 \begin{align*}
   g_{q,i,R}\otimes_r   g_{q,i,R} &= c_q^2 R^{d} \left\langle \int_{B_i}e_{xR}^{\otimes q} h(x) \, dx, \int_{B_i}e_{yR}^{\otimes q} h(y) \, dy \right\rangle_{\HH^{\otimes r}} \\
   &= c_q^2 R^{d}   \int_{B_i\times B_i}\big\langle e_{xR}^{\otimes q} ,   e_{yR}^{\otimes q}  \big\rangle_{\HH^{\otimes r}}        h(x)   h(y) \, dxdy \quad\text{by Fubini's theorem} \\
   &=c_q^2 R^{d}   \int_{B_i\times B_i}    \rho(xR-yR)^r  e_{xR}^{\otimes q-r} \otimes   e_{yR}^{\otimes q-r}       h(x)   h(y) \, dxdy
 \end{align*}
 and therefore,
\begin{align*}
& \big\| g_{q, i,R}\otimes_r g_{q, i,R} \big\| _{\HH^{\otimes 2q-2r}}^2 = c_q^4 R^{2d} \int_{B_i^4} \rho(Rx_1-Rx_2)^r  \rho(Rx_3-Rx_4)^r   \rho(Rx_1-Rx_3)^{q-r}   \\
& \qquad \qquad\qquad\qquad\qquad\qquad\qquad\qquad \qquad  \times   \rho(Rx_2-Rx_4)^{q-r} \prod _{i=1}^4h(x_i) \, d\pmb{x}   \\
 & = \frac {c_q^4} {R^{2d}}  \int_{ (RB_i)^4} \rho(x_1-x_2)^r  \rho(x_3-x_4)^r   \rho(x_1-x_3)^{q-r}    \rho(x_2-x_4)^{q-r}  
 \prod _{i=1}^4h(x_i/R) \, d\pmb{x} ,
 \end{align*}
 where $\pmb{x} =(x_1,x_2,x_3,x_4)$.
In view of the elementary inequality $a^r b^{q-r} \leq a^q + b^q$ for any $a, b\in\R_+$, we can write
\begin{align*}
&\quad \big\| g_{q, i,R}\otimes_r g_{q, i,R} \big\| _{\HH^{\otimes 2q-2r}}^2 \leq \frac{c_q^4}{R^{2d}} \int_{ (RB_i)^4}  \Big(  \big\vert  \rho(x_1-x_2) \big\vert^q  +\big\vert \rho(x_1-x_3)\big\vert^q \Big)  \\
&\qquad\qquad\qquad\qquad\qquad \qquad \times  \vert \rho(x_3-x_4)\vert^r   \vert   \rho(x_2-x_4)\vert^{q-r}   \prod _{i=1}^4 |h(x_i/R) | d\pmb{x}    \,.
\end{align*}
 Our goal is to show 
\[
\lim_{R\to+\infty} \sum_{r=1}^{q-1} \big\| g_{q, i,R}\otimes_r g_{q, i,R} \big\| _{\HH^{\otimes 2q-2r}}^2 = 0.
\]
Then by symmetry, it is enough to show that for each $r\in\{1, \dots, q-1\}$,
\begin{align*} 
\mathcal{K}_R= \frac{1}{ R^{2d}} \int_{ (RB_i)^4}    \vert  \rho(x_1-x_2) \vert^q     \vert \rho(x_3-x_4)\vert^r    \vert   \rho(x_2-x_4)\vert^{q-r} \prod _{i=1}^4|h(x_i/R) | d\pmb{x}\xrightarrow{R\to\infty} 0.
\end{align*}
Recall that $B_i$ is bounded, so we can assume $B_i\subset [-\ell, \ell]^d$ for some $\ell > 0$.   Taking  the continuity of $h$ into account  yields
\begin{align*}
\mathcal{K}_R&\leq  \frac{\big\| h\mathbf{1}_{[-\ell,\ell]^d}  \big\|_\infty}{ R^{2d}}\int_{ [-\ell R, \ell R]^{4d}}   \vert  \rho(x_1-x_2) \vert^q     \vert \rho(x_3-x_4)\vert^r    \vert   \rho(x_2-x_4)\vert^{q-r} 
 d\pmb{x} \\
&\leq    \frac{C}{R^{2d}} \left( \int_{\R^d} \vert \rho(z)\vert^q dz\right) \int_{ [-\ell R, \ell R]^{3d}}      \vert \rho(x_3-x_4)\vert^r   \vert   \rho(x_2-x_4)\vert^{q-r}   \, dx_2dx_3dx_4 \\
&\leq   \frac{C}{R^{d}} \left( \int_{\R^d} \vert \rho(z)\vert^m dz\right)   \left(\int_{[-2\ell R, 2\ell R]^{d}}      \vert \rho(x)\vert^r  dx\right)   \left( \int_{[-2\ell R, 2\ell R]^{d}} \vert   \rho(y)\vert^{q-r}  \, dy\right) .
\end{align*}

It suffices to show  that for each $r=1,\ldots, q-1$,  
\begin{align}\label{need001}
\frac{1}{R^{d(1 - rq^{-1})} }\int_{[-\ell R, \ell R  ]^d}      \vert \rho(x)\vert^r  dx \xrightarrow{R\to+\infty} 0 \,.
\end{align}
One can establish the above limit as follows.  Fix $\delta\in(0,1)$, we first decompose the above integral into two parts:  With $E(R) = [-\ell R, \ell R]^{d}$,
 \begin{align}\label{DEC}
 \frac{1}{R^{d(1 - rq^{-1})} }\int_{E(R)}      \vert \rho(x)\vert^r  dx  =\frac{\displaystyle \int_{E(\delta R)}      \vert \rho(x)\vert^r  dx}{R^{d(1 - rq^{-1})} }      + \frac{\displaystyle \int_{E(R) \setminus E(\delta R)}      \vert \rho(x)\vert^r  dx }{R^{d(1 - rq^{-1})} }  .
 \end{align}
 By H\"older's inequality, we have
 \[
  \int_{E(\delta R)}      \vert \rho(x)\vert^r  dx \leq  \left( \int_{\R^d}      \vert \rho(x)\vert^q  dx\right)^{r/q} ( 2\delta \ell R)^{d(1-rq^{-1})}
 \]
 and 
 \begin{align*}
  \int_{E(R) \setminus E(\delta R)}      \vert \rho(x)\vert^r  dx  &\leq  \left( \int_{E(R) \setminus E(\delta R)}      \vert \rho(x)\vert^q  dx \right)^{r/q} 
 2^{d(1-rq^{-1})}  \big[ (\ell R)^d -  (\delta\ell R)^d  \big]^{1-rq^{-1}} \\
  &\leq \left( \int_{\R^d} \mathbf{1}_{\{ \| x\| \geq \delta \ell R  \}}      \vert \rho(x)\vert^q  dx \right)^{r/q} ( 2\ell R)^{d(1 - rq^{-1})}.
  \end{align*}
Therefore, it is clear that due to $\rho\in L^q(\R^d, dx)$, for any fixed $\delta\in(0,1)$, the second term in \eqref{DEC} goes to zero, as $R\to+\infty$;  and the first term in \eqref{DEC} can be made arbitrarily small by choosing sufficiently small $\delta$. This completes our verification of condition (c) from Proposition \ref{CCLT}.

\medskip

(d) For each  $i\in\{1, \ldots, n\}$, we can see from the computations from step (a) that
\begin{align*}
 \sum_{q\geq N+1} q! \| g_{q, i,R} \|^2 _{\HH^{\otimes q}}  &=  \sum_{q\geq N+1}  c_q^2 q!  \int_{\{y\in B_i, y-zR^{-1} \in B_j \}} \rho(z)^q  h\big( y-zR^{-1}\big) h(y) \, dydz \\
 &\leq \left(\int_{\R^d}{\bf 1}_{B_i}(x) dx\right)   \left(\sup_{z\in B_i} |h(z)| ^2\right)  \sum_{q\geq N+1}  c_q^2 q!  \int_{\R^d} \vert \rho(z) \vert^m \, dz \,,
\end{align*}
which converges to zero (uniformly in $R$), as $N$ goes to infinity. 

Therefore, the limit in \eqref{FDD} is proved. In particular, \eqref{CLT1} is established.   \hfill $\square$

  \begin{remark}\label{rem33} If we only assume that $h:\R^d\to\R$ is continuous except at finitely many points, we can still obtain \eqref{CLT1}. This observation will be helpful in the proof of Theorem \ref{KU}.

  \end{remark}

    \subsection{Tightness} This part is split into two portions, dealing with proofs of \eqref{FCLT1} and \eqref{FCLT20} respectively.

   \begin{proof}[Proof of \eqref{FCLT1}] 
   For each $t\ge 0$, we   recall that $B_t = \{ x\in\R^d: \| x \| \leq t \}$ and   put
   \[
   X_R(t)=R^{d/2} \int_{B_t} \Phi(W_{xR}) h(x) dx.
   \]
Clearly  $X_R$ is a random variable with values in $C(\R_+)$. We know from Billingsley's book \cite{B99} that in order to have the tightness of $\{X_R, R>0\}$, it is sufficient to prove the following moment estimate: There exists some constant $C_T > 0$ such that for any $0 < s < t \leq T$,
        \begin{align}\label{MES}
    \big\| X_R(t) - X_R(s) \big\|_{L^p(\Omega)} \leq C_T \sqrt{ t -s   }\,,
    \end{align}
    where $p>2$ is  the fixed index  in the statement of Theorem \ref{Lawrence}. To simplify the  presentation, we assume that $T=1$.

Using the notation from Section \ref{PRE}, we first write $\Phi(W_{xR}) = \delta^m\Big( \Phi_m(W_{xR}) e_{xR}^{\otimes m} \Big)$. Then  for any $0 < s < t \leq 1$,
 \begin{align*}
  \big\| X_R(t) - X_R(s) \big\| _{L^p(\Omega)}    &=  R^{ d/2} \left\| \int_{B_{t} \setminus B_s } \Phi(W_{xR}) h(x) dx\right\| _{L^p(\Omega)} \\
 & = R^{ d/2} \left\| \int_{B_{t} \setminus B_s } \delta^m\Big(\Phi_m( W_{xR}) e_{xR}^{\otimes m}\Big) h(x) dx \right\| _{L^p(\Omega)}   \\
& = R^{ d/2} \left\| \delta^m\left(  \int_{B_{t} \setminus B_s }  \Phi_m( W_{xR}) e_{xR}^{\otimes m} ~ h(x)dx \right)\right\| _{L^p(\Omega)} \\
&= : \big\| \delta^m(v_R) \big\|   _{L^p(\Omega)}    \,,
   \end{align*}
with $v_R = R^{ d/2}  \int_{B_{t} \setminus B_s }  \Phi_m( W_{xR}) e_{xR}^{\otimes m} ~ h(x)dx $. Now we apply the Meyer's inequality (see \cite[Proposition 1.5.4]{Nualart06}), to get
\begin{align*}
  \big\| \delta^m(v_R) \big\|   _{L^p(\Omega)}  & \leq C \sum_{k=0}^m \big\| D^k v_R \big\| _{L^p(\Omega; \HH^{\otimes k+m})} \qquad\text{see also \cite[(2.8)]{NN18}} \\
 & \leq C \sum_{k=0}^m \left\|R^{ d/2}      \int_{B_{t} \setminus B_s } D^k \Big( \Phi_m( W_{xR}) e_{xR}^{\otimes m}\Big) ~ h(x)dx    \right\| _{L^p(\Omega; \HH^{\otimes (m+k)})} \,.
  \end{align*}
Keeping in mind the fact  that $h\in C(\R^d)$, we have 
   \begin{align}
   &  \left\|R^{ d/2}      \int_{B_{t} \setminus B_s } D^k \Big( \Phi_m( W_{xR}) e_{xR}^{\otimes m}\Big) ~ h(x)dx    \right\| _{L^p(\Omega; \HH^{\otimes (m+k)})}^2  \label{Quan01} \\
   &\leq   \Bigg\|  R^d    \int_{(B_{t} \setminus B_s )^2} \Big\langle  D^k\big(\Phi_m( W_{xR})\big), D^k\big(\Phi_m( W_{yR}) \big)  \Big\rangle_{\HH^{\otimes k}} \rho(xR-yR)^{m}    h(x)h(y)dxdy  \Bigg\| _{L^{\frac{p}{2}}(\Omega)} \notag \\
   &\leq  C  R^d  \int_{(B_{t} \setminus B_s )^2} \Big\|    \Big\langle  D^k\big(\Phi_m( W_{xR})\big), D^k\big(\Phi_m( W_{yR}) \big)   \Big\rangle_{\HH^{\otimes k}}    \Big\| _{L^{\frac{p}{2}}(\Omega)   } \vert \rho(xR-yR)\vert^{m} dxdy,\notag
       \end{align}
where we also applied Minkowski's inequality in the last inequality. Therefore, Cauchy-Schwarz inequality and  property (\pmb{B}) from Section \ref{PRE} imply that the quantity in \eqref{Quan01} is bounded by 
\begin{align*}
C  R^d  \int_{(B_{t} \setminus B_s )^2} \vert \rho(xR-yR)\vert^{m} dxdy  \leq C(t^d - s^d) \int_{\R^d} \vert \rho(z)\vert^m\, dz\,.
\end{align*}
 It follows that
 $
 \big\| X_R(t) - X_R(s) \big\| _{L^p(\Omega)}  \leq C \sqrt{t^d - s^d} \leq C \sqrt{t-s} \,.
  $
 \end{proof}

  \medskip

Now we show the weak convergence described in \eqref{FCLT20}.

   \begin{proof}[Proof of \eqref{FCLT20} ]  To simplify the notation, we assume $\ell=1$.  For $\pmb{z}\in  [0,1]^d$, we put
  \[
  Y_R(\pmb{z}) = R^{d/2} \int_{[0,\pmb{z}]}  \Phi(W_{xR}) h(x) dx
  \]
  and in what follows, we will  focus on establishing the tightness of $\{Y_R, R>0\}$ by proving the following estimate
  \begin{align}\label{es1}
 \| Y_R(\pmb{z})  - Y_R(\pmb{y})   \|_{L^p(\Omega)} \leq C \|  \pmb{z} - \pmb{y} \|^{1/2} \quad\text{for any $\pmb{y}, \pmb{z}\in [0,1]^d$,}
  \end{align}
   here $\|\cdot \|$ denotes the Euclidean norm and  $p > 2d$.   
    We write 
  \begin{align*}
  Y_R(\pmb{z})  - Y_R(\pmb{y}) &= R^{d/2}  \int_{ [0, \pmb{z} ] \setminus [0, \pmb{y} ]} g(W_{xR}) h(x) dx -  R^{d/2} \int_{ [0, \pmb{y} ] \setminus [0, \pmb{z} ]} g(W_{xR}) h(x) dx \\
  &= : A_1 - A_2 \,.
  \end{align*}
  Following the same arguments as in the proof of \eqref{FCLT1}, we have
  \begin{align*}
  \| A_1 \| _{L^p(\Omega)}^2& \leq C   R^d \int_{ [0, \pmb{z} ] \setminus [0, \pmb{y} ]} \int_{ [0, \pmb{z} ] \setminus [0, \pmb{y} ]}  \vert  \rho(Rx-Ry)\vert^{m}dxdy   \\
 &  \leq C \left( \int_{\R^d} \vert \rho(x) \vert^m dx\right) \max_{j=1}^d \vert y_j - z_j\vert  \leq C  \| \pmb{z} - \pmb{y} \|  .
   \end{align*}
 The same arguments yields the estimate
 $
  \| A_2 \| _{L^p(\Omega)} \leq C \| \pmb{y}- \pmb{z}\|^{1/2}$, so that \eqref{es1} holds true.  
 \end{proof}

\bigskip

\subsection{Proof of Theorem \ref{KU} } Put $q(x) = a(x)^{-1} - (1/a^\ast) = \Phi(W_x)$ and recall that  the solution to \eqref{EE} is given by 
\[
 u_\e(x) = c_\e(\omega) \int_0^x \frac{1}{a(y/\varepsilon)} \, dy - \int_0^x \frac{F(y)}{a(y/\varepsilon)} \, dy,
\] where
 $F(x) := \int_0^x f(y) \, dy$
 and 
 $$c_\e(\omega) : = \left( \, b + \int_0^1 \frac{F(y)}{a(y/\varepsilon)} \, dy \,\right) \left(  \int_0^1 \frac{1}{a(y/\varepsilon)} \, dy  \right)^{-1}.$$
 Note that for any $h\in C([0,1])$ and each $v\in(0,1]$, we obtain, by using the Hermite expansion, that
 \begin{align}
&\quad \left\| \int_0^v  q(y/\e)h(y) \, dy \right\|_{L^2(\Omega)}^2 = \sum_{q\geq m} c_q^2 q! \int_0^v\int_0^v \rho\big(\frac{y-x}{\e}\big)^q h(x) h(y) dxdy \notag \\
&\leq  \| h\|_\infty^2 \sum_{q\geq m} c_q^2 q! \int_0^v\int_0^v  \left| \rho\big(\frac{y-x}{\e}\big) \right|^m   dxdy \leq   \| h\|_\infty^2 \left(\sum_{q\geq m} c_q^2 q! \int_{\R} \vert\rho(z )\vert^m   dz\right)\e  \notag.
 \end{align}
 That is,
 \begin{align}
\left\| \int_0^v  q(y/\e)h(y) \, dy \right\|_{L^2(\Omega)}^2   \leq    \| h\|_\infty^2 \left(\sum_{q\geq m} c_q^2 q! \int_{\R} \vert\rho(z )\vert^m   dz\right)\e  \,. \label{need009}
 \end{align}
It follows that  
\begin{center}
$ {\displaystyle \int_0^v  \frac{1}{a(y/\e)} h(y) \, dy}$ converges in $L^2(\Omega)$ to ${\displaystyle \frac{1}{a^\ast} \int_0^v h(y) dy}$, as $\e\downarrow 0$. 
\end{center}
 In particular, the random vector
\[
J_\e(x) := \left(  \int_0^x \frac{1}{a(y/\varepsilon)} \, dy ,  \int_0^x \frac{F(y)}{a(y/\varepsilon)} \, dy , \int_0^1 \frac{F(y)}{a(y/\varepsilon)} \, dy,  \int_0^1 \frac{1}{a(y/\varepsilon)} \, dy \right)
\]
converges in $L^2(\Omega; \R^4)$ to 
\[
J(x) :=\left(  \frac{x}{a^\ast} ,  \int_0^x \frac{F(y)}{a^\ast} \, dy , \int_0^1 \frac{F(y)}{a^\ast} \, dy,   \frac{1}{a^\ast}  \right)\,.
\]
Put $M(z_1, z_2, z_3, z_4) = (b+z_3)z_1z_4^{-1} - z_2 $, then it follows from continuous mapping theorem  that   $u_{\e}(x) = M( J_\e(x))\to M( J(x)) =\bar{u}(x)$ in probability, as $\e\downarrow 0$, where 
 \[
 \bar{u}(x) = c^\ast \frac{x}{a^\ast} - \int_0^x \frac{F(y)}{a^\ast} dy \quad\text{with}\quad c^\ast= ba^\ast + \int_0^1 F(y) dy \,.
 \]
It is easy to see that $\bar{u}$ solve  equation \eqref{EEh}, so  part (1) of Theorem \ref{KU} is established.   

Following the decomposition given in \cite[pages 1082-1085]{GB12},    we rewrite the rescaled corrector as follows:
\begin{align}\label{deomp_0}
  \frac{u_\e(x)-\bar{u}(x)}{ \sqrt{\e} } = \mathcal{U}_\e(x) +        r_\e(x) +   \frac{ \rho_\e(x)  }{ \sqrt{\e} }\, ,
\end{align}
where ${\displaystyle  
\mathcal{U}_\e(x) := \frac{1}{\sqrt{\e}} \int_0^1 F(x, y) q(y/\e)  \, dy}$, ${\displaystyle  r_\e(x) := \frac{c_\e-c^\ast}{\sqrt{\e}}\int_0^x q(y/\varepsilon) \, dy}$ and
\begin{align*}
  \rho_\e(x) & := \frac{x}{\displaystyle\int_0^1 \dfrac{1}{a(y/\varepsilon)} \, dy}   \left[ c^\ast\left(\int_0^1 q(y/\varepsilon)  \, dy\right)^2   - \left(\int_0^1F(y) q(y/\varepsilon)  \, dy\right)    \int_0^1 q(y/\varepsilon)  \, dy \right],
\end{align*}
with $F(x,y) = \big( c^\ast - F(y) \big)\1_{[0,x]}(y) + x \big( F(y) - c^\ast\big)\1_{[0,1]}(y)$.  Therefore, it follows from Theorem \ref{Lawrence} and the observation in Remark \ref{rem33} that  $\mathcal{U}_{\e}(x)$ converges to a centered Gaussian distribution with variance $\mu^2\int_0^1 F(x,y)^2 dy$.
Let us show that the terms  $r_{\e}(x)$ and $\rho_{\e}(x)$ do not contribute to the limit.

\medskip

\noindent{(i) \it Estimation of $r_{\e}(x)$}: We know   that $\e^{-1/2} \int_0^x q(y/\e) \, dy$ converges in law to a Gaussian random variable and  $c_{\e} - c^\ast$ converges in probability to zero, as $\e\downarrow 0$. 
It follows   that $r_{\e}(x)$ converges in probability to zero, as $\e\downarrow 0$. Moreover, under the additional assumption that $\Phi\in L^p(\R, e^{-x^2/2}dx)$ with $p>2$, we can  apply \eqref{FCLT20} with $d=1$ and $R= 1/\e$ and conclude that, as $\e\downarrow 0$,
\[
\left\{    \frac{1}{\sqrt{\e}}\int_0^x q(y/\e) \, dy\,, x\in[0,1] \right\} \quad\text{converges in law to a Gaussian process.}
\]
Thus,  the process $\{r_{\e}(x),  x\in[0,1] \}$ converges in law, hence also in probability, to the zero process.

\medskip

\noindent{(ii) \it Estimation of $\rho_{\e}(x)$}: Similarly,  
\begin{align*}
 \frac{ \rho_{\e}(x)  }{ \sqrt{\e} } =&  c^\ast \sqrt{\e} \left(\int_0^1 \dfrac{1}{a(y/\e)} \, dy\right)^{-1}  \left( \frac{1}{\sqrt{\e}}\int_0^1 q(y/\e)  \, dy\right)^2   x  \\
  &\quad   -  \left(   \int_0^1 q(y/\e)  \, dy \right)\left(  \int_0^1 \dfrac{1}{a(y/\e)} \, dy\right)^{-1}  \left(\frac{1}{\sqrt{\e}}\int_0^1F(y) q(y/\e)  \, dy\right)  x\,.  
\end{align*}
 It is clear that  both
  $$c^\ast \sqrt{\e} \left(\int_0^1 \dfrac{1}{a(y/\e)} \, dy\right)^{-1}   \quad{\rm  and}\quad    \left( \int_0^1 q(y/\e)  \, dy\right)\left( \int_0^1 \dfrac{1}{a(y/\e)} \, dy\right)^{-1}$$ 
   converge to zero in probability, while  both 
 $$ \left\{ \left( \frac{1}{\sqrt{\e}}\int_0^1 q(y/\e)  \, dy\right)^2 \cdot x \right\}_{x\in[0,1]} \quad \text{and}\quad  \left\{ \left(\frac{1}{\sqrt{\e}}\int_0^1F(y) q(y/\e)  \, dy\right) \cdot x   \right\}_{x\in[0,1]} $$
weakly converge to some processes in $C([0,1])$, as $\e\downarrow 0$.
This implies the process $\big\{ \e^{-1/2} \rho_{\e}(x)\,: x\in[0,1] \big\}$ converges in probability to  the zero process. Then it follows immediately that   
\begin{center}
  $ r_{\e}(x) +   \dfrac{ \rho_{\e}(x)  }{ \sqrt{\e} }$  converges in probability to the zero, for every $x\in[0,1]$;
\end{center}
and  under the additional assumption   $\Phi\in L^p(\R, e^{-x^2/2}dx)$,
\begin{center}
  $\left\{ r_{\e}(x) +   \dfrac{ \rho_{\e}(x)  }{ \sqrt{\e} }, x\in[0,1] \right\}$  converges in probability to the zero process.
\end{center}

\medskip

\noindent{(iii)\it Endgame:} In view of Slutsky's theorem, we have just established \eqref{goal1} and to reach \eqref{goal2},  it suffices to prove  as $\e\downarrow 0$,
$$
\mathcal{U}_{\e} \xrightarrow{\rm law} \left\{  \mu \int_{0}^1 F(x,y) dA_y\, x\in [0,1] \right\},
$$
 where  $A$  is  a standard Brownian motion on $[0,1]$.  Now we write for every $x\in[0,1]$,
\begin{align*}
\mathcal{U}_{\e}(x) &= \frac{1}{\sqrt{\e}} \int_0^x \big( c^\ast - F(y) \big) q(y/\e) dy + \frac{x}{\sqrt{\e}} \int_0^1 \big(F(y) - c^\ast \big) q(y/\e)dy \\
&=: \mathcal{V}_{1,\e}(x) + \mathcal{V}_{2,\e}(x) \,.
\end{align*}
Then applying \eqref{FCLT20}  again yields
\begin{align*}
\mathcal{V}_{1,\e} \xrightarrow[\rm law]{\e\downarrow 0} \mathcal{V}_1 :=   \left\{  \mu \int_0^x \big( c^\ast - F(y) \big) dA_y, x\in [0,1] \right\}   \,.
\end{align*}
  Note that we can write $\mathcal{V}_{2,\e}(x) = x \kappa_\e$ with $\kappa_\e = \e^{-1/2} \int_0^1 \big(F(y) - c^\ast \big) q(y/\e)dy $ bounded in $L^2(\Omega)$ in view of \eqref{need009}.  It is also clear that  
 as $\e\downarrow 0$, $\kappa_\e$ converges in law to $  \mu  \int_0^1  \big( F(y) - c^\ast \big) dA_y$. As a consequence, the \emph{f.d.d.} convergence of $\mathcal{V}_{2,\e}$ is trivial and the tightness follows from the fact that
$$\E\big[ (\mathcal{V}_{2,\e}(x)  - \mathcal{V}_{2,\e}(y) \big)^2 \big] = |x-y|^2 \E\big[ \kappa_\e ^2 \big] \leq C |x-y|^2 \quad\text{by \eqref{need009}}.
$$
Thus,
\[
\mathcal{V}_{2,\e} \xrightarrow{\rm law}  \mathcal{V}_2 :=  \left\{   \mu x \int_0^1  \big( F(y) - c^\ast \big) dA_y \,,\, x\in[0,1] \right\} .
\]

It follows that the sequence $(\mathcal{V}_{1,\e}, \mathcal{V}_{2,\e}   )$ is tight, and so is $\mathcal{V}_{1,\e}+ \mathcal{V}_{2,\e} $. That is, $\mathcal{U}_{\e}$ is tight.  Now consider $\lambda_k\in\R$ and $x_k\in[0,1]$ for $k\in\{1,\ldots, \ell\}$ and any $\ell\geq 1$. We have
\[
\sum_{k=1}^\ell  \lambda_k\mathcal{U}_{\e}(x_k) =  \frac{1}{\sqrt{\e}} \int_0^1 \sum_{k=1}^\ell  \lambda_kF(x_k, y) q(y/\e)  \, dy\xrightarrow[\rm law]{\e\downarrow 0}\sum_{k=1}^\ell  \lambda_k ~\sigma \int_0^1  F(x_k, y)   \, dA_y\,.
\]
This proves the convergence of  the finite-dimensional distributions for $\mathcal{U}_{\e}$ and conclude our proof with the above tightness of $\big\{\mathcal{U}_{\e}, \e > 0 \big\}$.   \hfill $\square$

\section{Proof of Theorem \ref{Kansas}} \label{SEC4}

The proof follows similar arguments as in the proof of Theorem \ref{Lawrence}. Here we first sketch the proof of \eqref{CLT2}.  For any $B\in \mathcal{B}_b$, we first rewrite using Hermite expansions
\begin{align*}
\widehat{G}_R(B):&= R^{\alpha/2} \int_{B} \Phi(W_{xR}) \nu(dx) =   R^{\alpha/2} \int_{B} \sum_{q\geq m}c_q H_q(W_{xR}) \nu(dx) \\
&   =     \sum_{q\geq m} \delta^q\left(  c_q R^{\alpha/2} \int_{B} e_{xR}^{\otimes q} \nu(dx)  \right).
\end{align*}
By the orthogonality of Hermite polynomials, we have 
\begin{align*}
\E\big[ \widehat{G}_R(B)^2 \big] &=  R^\alpha \sum_{q\geq m} c_q^2 q! \int_{B^2} \rho(xR-yR)^q \nu(dx) \nu(dy) \\
&=  R^{-\alpha} \sum_{q\geq m} c_q^2 q! \int_{(RB)^2} \rho(x-y)^q \nu(dx) \nu(dy)\,,
\end{align*}
where we used the $\alpha$-homogeneity and made a change of variable in the last equality:  $(xR, yR)\to (x,y)$.     Making another change of variable $(x=x, z=x-y)$ yields
\begin{align*}
\int_{(RB)^2} \rho(x-y)^q \nu(dx) \nu(dy)& = \int_{\R^d} \rho(z)^q \nu\big[ (RB)\cap (z+RB)  \big] \, \nu(dz) \\
&= R^\alpha \nu(B) \int_{\R^d} \rho(z)^q  \frac{\nu\big[ (RB)\cap (z+ RB)  \big] }{\nu(RB)}\, \nu(dz) \,.
\end{align*}
In view of the $\alpha$-homogeneity, the quantity $\nu [ (RB)\cap (z+ RB)   ] / \nu(RB)$ converges to $1$ as $R\to+\infty$, for each $z\in\R^d$. Indeed, given $z\in\R$, we can write
\[
\frac {\nu [ (RB)\cap (z+ RB)   ] }{\nu(RB)}= \frac{ \nu [ B\cap (R^{-1}z+ B)   ] }{ \nu(B)}.
\]
 By the dominated convergence theorem, our assumptions ensure that
\begin{align*}
R^\alpha \E\big[\widehat{G}_R(B)^2\big] &=    \nu(B) \sum_{q\geq m} c_q^2 q!     \int_{\R^d} \rho(z)^q  \frac{\nu\big[ (RB)\cap  (z+ RB)  \big] }{\nu(RB)}\, \nu(dz) \\
&\xrightarrow{R\to+\infty}  \nu(B) \sum_{q\geq m} c_q^2 q!     \int_{\R^d} \rho(z)^q  \, \nu(dz) \,.
\end{align*}
This gives us the limiting variance. To   show the central convergence,  it is routine to verify the contraction conditions, which can be done in the same way as before. We omit the details here and point out that we   need to use the following limiting result instead of \eqref{need001}: for each $r\in\{ 1,\ldots, q-1\}$,  
\[
\frac{1}{R^{\alpha(1 - rq^{-1})} }\int_{RB}      \vert \rho(x)\vert^r  \nu(dx) \xrightarrow{R\to+\infty} 0 \,. 
\]
The above limit can be verified in the same way, by using H\"older's inequality and the fact that $\rho\in L^m(\R^d, d\nu)$.  In this way, we can obtain the \emph{f.d.d.} convergence described in \eqref{CLT2}, and we leave this as an easy exercise for the interested readers.

\medskip

 In the following, we sketch the arguments for tightness.  For every $t\ge 0$, we put 
\[
\widehat{X}_R(t ) = R^{ \alpha/2} \int_{B_t} \Phi(W_{xR}) \nu(dx).
\]
In the sequel, we 
    show  the tightness for  $\{ \widehat{X}_R, R > 0 \}$. 
 
 For fixed $0< s  < t \leq 1$, we can obtain,  by similar arguments  as before, that
 \begin{align*}
&\quad \big\|  \widehat{X}_R(t) - \widehat{X}_R(s) \big\| _{L^p(\Omega)} =  R^{\alpha/2} \left\| \delta^m\left( \int_{B_t\setminus B_s} \Phi_m(W_{xR}) e^{\otimes m}_{xR} \nu(dx)  \right) \right\| _{L^p(\Omega)}\\
& \leq C    \left( R^{\alpha} \int_{\big(B_{t}\setminus B_{s}\big)^2}  \vert \rho(xR-yR)\vert^{m} \nu(dx) \nu(dy)  \right)^{1/2}  \\
&\leq C   \left(   R^{-\alpha}  \nu\Big(  B_{tR}\setminus B_{sR} \Big)    \int_{\R^d }  \vert \rho(z)\vert^{m}  \nu(dz) \right)^{1/2}  \leq C  \big( t^\alpha - s^\alpha \big)^{1/2} \,;
 \end{align*}
 see proof of \eqref{FCLT1}. Note that for any $a,b\in\R_+$ and any $\beta\in(0,1]$, it holds that $(a+b)^\beta \leq a^\beta + b^\beta$;  for any $a,b\in[0,1]$ and $\beta\in (1,+\infty)$, there exists a constant $C_\beta$ that only depends on $\beta$ such that $\vert a^\beta - b^\beta \vert \leq C_\beta \vert a-b\vert$. 
This gives us 
 \[
\big\|  \widehat{X}_R(t) - \widehat{X}_R(s) \big\| _{L^p(\Omega)} \leq     C (t-s)^{( \alpha\wedge 1) /2} \,.
 \]
Since $\alpha p > 2$, we can deduce the tightness of $\{ \widehat{X}_R, R > 0 \}$. \hfill $\square$

 \medskip
 \noindent
 {\bf  Acknowledgement:} We would like to thank two anonymous referees for their valuable suggestions  that help us to improve the presentation of our results.

\end{document}